\newcommand{\nc}{\newcommand}
\nc{\rnc}{\renewcommand}

\newtheorem{proposition}{Proposition}

\newcommand{\pf}{\noindent{\bf Proof}\,\,  }
\newcommand{\epf}{\hfill$\Box$\bigskip}

\nc{\pfnb}{\pf} \nc{\epfnb}{\bigskip}

\newtheorem{thm}[equation]{Theorem}
\newtheorem{lemma}[equation]{Lemma}

\rnc{\th}{\theta} \rnc{\H}{\mathscr H} \nc{\X}{\mathscr X}
\nc{\Y}{\mathscr Y} \nc{\XR}{\pres{\X}{\R\,}}
\nc{\rev}{\text{rev}} \nc{\suvertex}[1]{\pscircle*(#1,2){1.5 pt}}
\nc{\slvertex}[1]{\pscircle*(#1,0){1.5 pt}}
\nc{\AND}{\qquad\text{and}\qquad}

\nc{\uv}[1]{\pscircle*(#1,2){2pt}}
\nc{\lv}[1]{\pscircle*(#1,0){2pt}}
\nc{\stline}[2]{\psline(#1,2)(#2,0)}
\nc{\guv}[1]{\pscircle*[linecolor=lightgray](#1,2){2pt}}
\nc{\glv}[1]{\pscircle*[linecolor=lightgray](#1,0){2pt}}
\nc{\gstline}[2]{\psline[linecolor=lightgray](#1,2)(#2,0)}
\nc{\ul}[2]{\psline(#1,2)(#2,2)}
\rnc{\ll}[2]{\psline(#1,0)(#2,0)}
\nc{\uc}{\pscurve(0,2)(1,1.6)(2,2)}
\nc{\lc}{\pscurve(0,0)(1,0.4)(2,0)}
\nc{\gul}[2]{\psline[linecolor=lightgray](#1,2)(#2,2)}
\nc{\gll}[2]{\psline[linecolor=lightgray](#1,0)(#2,0)}
\nc{\guc}{\pscurve[linecolor=lightgray](0,2)(1,1.6)(2,2)}
\nc{\glc}{\pscurve[linecolor=lightgray](0,0)(1,0.4)(2,0)}
\nc{\uvtx}[1]{\pscircle*(#1,4){2pt}}
\nc{\uvertex}[1]{\pscircle*(#1,2){2pt}}
\nc{\lvertex}[1]{\pscircle*(#1,0){2pt}}
\nc{\stlinebd}[2]{\psline[border=0.5 mm](#1,2)(#2,0)}

\nc{\guvertex}{\guv} \nc{\glvertex}{\glv}

\newcommand{\larr}{\mathrel{\hspace{-0.35 ex}>\hspace{-1.1ex}-}\hspace{-0.35 ex}}
\newcommand{\rarr}{\mathrel{\hspace{-0.35ex}-\hspace{-0.5ex}-\hspace{-2.3ex}>\hspace{-0.35 ex}}}
\newcommand{\lrarr}{\mathrel{\hspace{-0.35ex}>\hspace{-1.1ex}-\hspace{-0.5ex}-\hspace{-2.3ex}>\hspace{-0.35 ex}}}

\nc{\bit}{\vspace{-3 truemm}\begin{itemize}}
\nc{\eit}{\end{itemize}\vspace{-3 truemm}}

\nc{\arw}[1]{\psset{xunit=1cm, yunit=1cm}
 \begin{pspicture}(0,0)(0.52,0.2)
  \psline[linewidth=0.5pt]{#1}(0.05,0.1)(0.47,0.1) \end{pspicture}}

\documentclass{article}

\usepackage{amsmath,amssymb,latexsym,pstricks,mathrsfs,xypic}
\begin{document}

\title{A Presentation for the Dual Symmetric Inverse Monoid}
\author{David Easdown\\
{\footnotesize \emph{School of Mathematics and Statistics, University of Sydney, NSW 2006, Australia}}\\
{\footnotesize {\tt de\,@\,maths.usyd.edu.au} }\\~\\
James East\\
{\footnotesize \emph{Department of Mathematics, La Trobe
University, Victoria 3083, Australia}}\\
{\footnotesize {\tt james.east\,@\,latrobe.edu.au} }\\~\\
D.~G.~FitzGerald\\
{\footnotesize \emph{School of Mathematics and Physics, University of Tasmania, Private Bag 37, Hobart 7250, Australia}}\\
{\footnotesize {\tt d.fitzgerald\,@\,utas.edu.au} }}

\maketitle


\begin{abstract}
The dual symmetric inverse monoid $\mathscr{I}_n^*$ is the inverse monoid
of all isomorphisms between quotients of an $n$-set. We give a
monoid presentation of $\mathscr{I}_n^*$ and, along the way, establish
criteria for a monoid to be inverse when it is generated by
completely regular elements.
\end{abstract}



\section{Introduction}

Inverse monoids model the partial or local symmetries of structures, generalizing
the total symmetries modelled by groups. Key examples are the \emph{symmetric
inverse monoid} $\mathscr{I}_X$ on a set $X$ (consisting of all bijections between subsets
of $X$), and the \emph{dual symmetric inverse monoid} $\mathscr{I}_X^*$ on $X$ (consisting
of all isomorphisms between subobjects of $X$ in the category ${\bf
Set}^{\text{opp}}$), each with an appropriate multiplication. They share the
property that every inverse monoid may be faithfully represented in some $\mathscr{I}_X$
and some $\mathscr{I}_X^*$. The monoid $\mathscr{I}_X^*$ may be realized in many different ways; in
\cite{2}, it was described as consisting of bijections between quotient sets of
$X$, or \emph{block bijections} on $X$, which map the blocks of a ``domain''
equivalence (or partition) on $X$ bijectively to blocks of a ``range''
equivalence. These objects may also be regarded as special binary relations on $X$
called \emph{biequivalences}. The appropriate multiplication involves the join of
equivalences---details are found in \cite{2}, and an alternative description in
\cite[pp.~122--124]{4}.

\subsection{Finite dual symmetric inverse
monoids}\label{sect:fdsim}

In this paper we focus on finite $X$, and write $\mathbf{n}=\{  1,\dots n\}$ and
$\mathscr{I}_n^*=\mathscr{I}_{\mathbf{n}}^*$. In a graphical representation described in \cite{5}, the
elements of $\mathscr{I}_n^*$ are thought of as graphs on a vertex
set~$\{1,\ldots,n\}\cup\{1',\ldots,n'\}$ (consisting of two copies of $\mathbf{n}$) such
that each connected component has at least one dashed and one undashed vertex.
This representation is not unique---two graphs are regarded as equivalent if they
have the same connected components---but it facilitates visualization and is
intimately connected to the combinatorial structure. Conventionally, we draw the
graph of an element of $\mathscr{I}_n^*$ such that the vertices $1,\ldots,n$ are in a
horizontal row (increasing from left to right), with vertices $1',\ldots,n'$
vertically below. See Fig. 1 for the graph of a block bijection $\theta\in \mathscr{I}_{8}^*$ 
with domain $(1,2\,|\,3\,|\,4,6,7\,|\,5,8)$ and range
$(1\,|\,2,4\,|\,3\,|\,5,6,7,8)$. In an obvious notation, we also write $\textstyle{{\theta = \big( {1,2 \atop 2,4} \big|
{3 \atop 5,6,7,8} \big| {4,6,7 \atop 1} \big| {5,8 \atop 3} \big).}}$

\begin{figure}[ht]
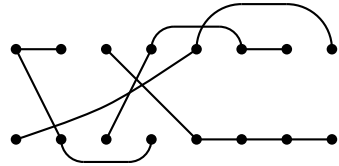

   \begin{center}
 \vspace{2 cm}
 \psset{unit= .6 cm}
 \psset{origin={-3.5,0}}
 \pscircle*(0,2){2pt}
 \pscircle*(1,2){2pt}
 \pscircle*(2,2){2pt}
 \pscircle*(3,2){2pt}
 \pscircle*(4,2){2pt}
 \pscircle*(5,2){2pt}
 \pscircle*(6,2){2pt}
 \pscircle*(7,2){2pt}
 \pscircle*(0,0){2pt}
 \pscircle*(1,0){2pt}
 \pscircle*(2,0){2pt}
 \pscircle*(3,0){2pt}
 \pscircle*(4,0){2pt}
 \pscircle*(5,0){2pt}
 \pscircle*(6,0){2pt}
 \pscircle*(7,0){2pt}
 \psline(0,2)(1,2)
 \psline(5,2)(6,2)
 \psline(4,0)(7,0)
 \psline(0,2)(1,0)
 \psline(2,2)(4,0)
 \psline(3,2)(2,0)
 \pscurve(4,2)(2,.75)(0,0)
 \psarc[linewidth=.8pt](3.5,2){.5}{90}{180}
 \psarc[linewidth=.8pt](4.5,2){.5}{0}{90}
 \psline[linewidth=.8pt](3.5,2.5)(4.5,2.5)
 \psarc[linewidth=.8pt](5,2){1}{90}{180}
 \psarc[linewidth=.8pt](6,2){1}{0}{90}
 \psline[linewidth=.8pt](5,3)(6,3)
 \psarc[linewidth=.8pt](1.5,0){.5}{180}{270}
 \psarc[linewidth=.8pt](2.5,0){.5}{270}{0}
 \psline[linewidth=.8pt](1.5,-0.5)(2.5,-0.5)
 \caption{A graphical representation of a block bijection $\theta \in \mathscr{I}_8^*$.}
 \label{picoftheta}
   \end{center}
 \end{figure}
To multiply two such diagrams, they are stacked
vertically, with the ``interior'' rows of vertices coinciding; then the connected
components of the resulting graph are constructed and the interior vertices are
ignored. See Fig. 2 for an example.
\begin{figure}[ht]
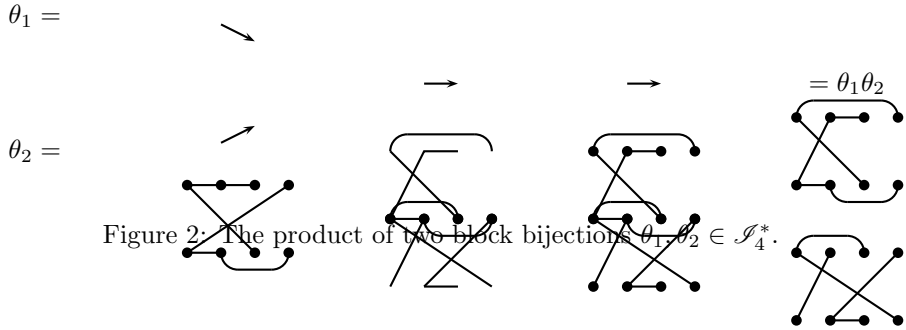

   \begin{center}
 \vspace{3.3 cm}
 \psset{xunit= .45 cm,yunit= .45 cm,runit=.45 cm}
 \psset{origin={-10.5,4}}
 \pscircle*(0,2){2pt}
 \pscircle*(1,2){2pt}
 \pscircle*(2,2){2pt}
 \pscircle*(3,2){2pt}
 \pscircle*(0,0){2pt}
 \pscircle*(1,0){2pt}
 \pscircle*(2,0){2pt}
 \pscircle*(3,0){2pt}
 \psline(1,2)(0,0)
 \psline(0,2)(3,0)
 \psline(3,2)(1,0)
 \psline(1,0)(2,0)
 \psarc[linewidth=.8pt](0.5,2){.5}{90}{180}
 \psarc[linewidth=.8pt](1.5,2){.5}{0}{90}
 \psline[linewidth=.8pt](0.5,2.5)(1.5,2.5)
 \psset{origin={-10.5,0}}
 \pscircle*(0,2){2pt}
 \pscircle*(1,2){2pt}
 \pscircle*(2,2){2pt}
 \pscircle*(3,2){2pt}
 \pscircle*(0,0){2pt}
 \pscircle*(1,0){2pt}
 \pscircle*(2,0){2pt}
 \pscircle*(3,0){2pt}
 \psline(2,2)(1,2)(0,0)(1,0)
 \psline(0,2)(2,0)
 \psarc[linewidth=.8pt](0.5,2){.5}{90}{180}
 \psarc[linewidth=.8pt](2.5,2){.5}{0}{90}
 \psline[linewidth=.8pt](0.5,2.5)(2.5,2.5)
 \psarc[linewidth=.8pt](1.5,0){.5}{180}{270}
 \psarc[linewidth=.8pt](2.5,0){.5}{270}{0}
 \psline[linewidth=.8pt](1.5,-0.5)(2.5,-0.5)
 \psset{origin={-4.5,3}}
 \pscircle*(0,2){2pt}
 \pscircle*(1,2){2pt}
 \pscircle*(2,2){2pt}
 \pscircle*(3,2){2pt}
 \pscircle*(0,0){2pt}
 \pscircle*(1,0){2pt}
 \pscircle*(2,0){2pt}
 \pscircle*(3,0){2pt}
 \psline(1,2)(0,0)
 \psline(0,2)(3,0)
 \psline(3,2)(1,0)
 \psline(1,0)(2,0)
 \psarc[linewidth=.8pt](0.5,2){.5}{90}{180}
 \psarc[linewidth=.8pt](1.5,2){.5}{0}{90}
 \psline[linewidth=.8pt](0.5,2.5)(1.5,2.5)
 \psset{origin={-4.5,1}}
 \pscircle*(0,2){2pt}
 \pscircle*(1,2){2pt}
 \pscircle*(2,2){2pt}
 \pscircle*(3,2){2pt}
 \pscircle*(0,0){2pt}
 \pscircle*(1,0){2pt}
 \pscircle*(2,0){2pt}
 \pscircle*(3,0){2pt}
 \psline(2,2)(1,2)(0,0)(1,0)
 \psline(0,2)(2,0)
 \psarc[linewidth=.8pt](0.5,2){.5}{90}{180}
 \psarc[linewidth=.8pt](2.5,2){.5}{0}{90}
 \psline[linewidth=.8pt](0.5,2.5)(2.5,2.5)
 \psarc[linewidth=.8pt](1.5,0){.5}{180}{270}
 \psarc[linewidth=.8pt](2.5,0){.5}{270}{0}
 \psline[linewidth=.8pt](1.5,-0.5)(2.5,-0.5)
 \psset{origin={1.5,3}}
 \pscircle*(0,2){2pt}
 \pscircle*(1,2){2pt}
 \pscircle*(2,2){2pt}
 \pscircle*(3,2){2pt}
 \psline(1,2)(0,0)
 \psline(0,2)(3,0)
 \psline(3,2)(1,0)
 \psline(1,0)(2,0)
 \psarc[linewidth=.8pt](0.5,2){.5}{90}{180}
 \psarc[linewidth=.8pt](1.5,2){.5}{0}{90}
 \psline[linewidth=.8pt](0.5,2.5)(1.5,2.5)
 \psset{origin={1.5,1}}
 \pscircle*(0,0){2pt}
 \pscircle*(1,0){2pt}
 \pscircle*(2,0){2pt}
 \pscircle*(3,0){2pt}
 \psline(2,2)(1,2)(0,0)(1,0)
 \psline(0,2)(2,0)
 \psarc[linewidth=.8pt](0.5,2){.5}{90}{180}
 \psarc[linewidth=.8pt](2.5,2){.5}{0}{90}
 \psline[linewidth=.8pt](0.5,2.5)(2.5,2.5)
 \psarc[linewidth=.8pt](1.5,0){.5}{180}{270}
 \psarc[linewidth=.8pt](2.5,0){.5}{270}{0}
 \psline[linewidth=.8pt](1.5,-0.5)(2.5,-0.5)
 \psset{origin={7.5,2}}
 \pscircle*(0,2){2pt}
 \pscircle*(1,2){2pt}
 \pscircle*(2,2){2pt}
 \pscircle*(3,2){2pt}
 \pscircle*(0,0){2pt}
 \pscircle*(1,0){2pt}
 \pscircle*(2,0){2pt}
 \pscircle*(3,0){2pt}
 \psline(2,2)(0,2)(2,0)
 \psline(3,2)(0,0)(1,0)
 \psarc[linewidth=.8pt](1.5,0){.5}{180}{270}
 \psarc[linewidth=.8pt](2.5,0){.5}{270}{0}
 \psline[linewidth=.8pt](1.5,-0.5)(2.5,-0.5)
 \psset{origin={0,0}}
 \psline{->}(-6.5,4.75)(-5.5,4.25)
 \psline{->}(-6.5,1.25)(-5.5,1.75)
 \psline{->}(-.5,3)(.5,3)
 \psline{->}(5.5,3)(6.5,3)
 \rput(-12,5){$\theta_1=$}
 \rput(-12,1){$\theta_2=$}
 \rput(12,3){$=\theta_1\theta_2$}
 \caption{The product of two block bijections $\theta_1,\theta_2\in\mathscr{I}_4^*$.}\label{prodininstar}
   \end{center}
 \end{figure}
\newline It is clear from its graphical representation that
$\mathscr{I}_n^*$ is a submonoid of the \emph{partition monoid}, though not one of the
submonoids discussed in \cite{3}. Maltcev \cite{5} shows that $\mathscr{I}_n^*$ with the
zero of the partition monoid adjoined is a maximal inverse subsemigroup of the
partition monoid, and gives a set of generators for $\mathscr{I}_n^*$. These generators are
completely regular; later in this paper, we present auxiliary results on the
generation of inverse semigroups by completely regular elements. Although these
results are of interest in their own right, our main goal is to obtain a
presentation, in terms of generators and relations, of~$\mathscr{I}_n^*$. Our method makes
use of known presentations of some special subsemigroups of $\mathscr{I}_n^*$. We now
describe these subsemigroups, postponing their presentations until a later
section.

The group of units of $\mathscr{I}_n^*$ is the symmetric group $\mathcal{S}_n$, while the
semilattice of idempotents is (isomorphic to) $\mathscr{E}_n$, the set of all equivalences
on $\mathbf{n}$, with multiplication being join of equivalences. Another subsemigroup
consists of those block bijections which are induced by permutations of $\mathbf{n}$
acting on the equivalence relations; this is the \emph{factorizable part}
of~$\mathscr{I}_n^*$, which we denote by $\mathscr{F}_n$, and which is equal to the set product
$\mathscr{E}_n\mathcal{S}_n=\mathcal{S}_n\mathscr{E}_n$. In~\cite{2} these elements were called \emph{uniform}, and
in \cite{5} \emph{type-preserving}, since they have the characteristic property
that corresponding blocks are of equal cardinality. We will also refer to the
\emph{local submonoid}~$\varepsilon\mathscr{I}_X^*\varepsilon$ of $\mathscr{I}_X^*$ determined by a non-identity
idempotent $\varepsilon$. This subsemigroup consists of all $\beta\in\mathscr{I}_X^*$ for which $\varepsilon$
is a (left and right) identity. Recalling that the idempotent~$\varepsilon$ is an
equivalence on $X$, it is easy to see that there is a natural
isomorphism~$\varepsilon\mathscr{I}_X^*\varepsilon\to\mathscr{I}_{X/\varepsilon}^*$. As an example which we make use of
later, when $X=\mathbf{n}$ and~$\varepsilon=(1,2\,|\,3\,|\,\cdots\,|\,n)$, we obtain an
isomorphism $\Upsilon:\varepsilon\mathscr{I}_n^*\varepsilon\to\mathscr{I}_{n-1}^*$. Diagrammatically, we obtain a
graph of $\beta\Upsilon\in\mathscr{I}_{n-1}^*$ from a graph of $\beta\in\varepsilon\mathscr{I}_n^*\varepsilon$ by
identifying vertices $1\equiv2$ and $1'\equiv2'$, relabelling the vertices, and
adjusting the edges accordingly; an example is given in Fig. 3.
\begin{figure}[ht]
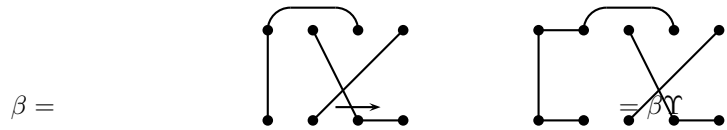

   \begin{center}
 \vspace{2 cm}
 \psset{xunit= .6 cm,yunit= .6 cm,runit=.6 cm}
 \psset{origin={-5,0}}
 \pscircle*(0,2){2pt}
 \pscircle*(1,2){2pt}
 \pscircle*(2,2){2pt}
 \pscircle*(3,2){2pt}
 \pscircle*(0,0){2pt}
 \pscircle*(1,0){2pt}
 \pscircle*(2,0){2pt}
 \pscircle*(3,0){2pt}
 \uvertex{-1}
 \lvertex{-1}
 \psline(-1,2)(0,2)
 \psline(1,2)(2,0)(3,0)
 \psline(-1,0)(0,0)
 \psline(3,2)(1,0)
 \stline{-1}{-1}
 \psarc[linewidth=.8pt](0.5,2){.5}{90}{180}
 \psarc[linewidth=.8pt](1.5,2){.5}{0}{90}
 \psline[linewidth=.8pt](0.5,2.5)(1.5,2.5)
 \psset{origin={2,0}}
 \pscircle*(0,2){2pt}
 \pscircle*(1,2){2pt}
 \pscircle*(2,2){2pt}
 \pscircle*(3,2){2pt}
 \pscircle*(0,0){2pt}
 \pscircle*(1,0){2pt}
 \pscircle*(2,0){2pt}
 \pscircle*(3,0){2pt}
 \psline(1,2)(2,0)(3,0)
 \psline(0,0)(0,2)
 \psline(3,2)(1,0)
 \psarc[linewidth=.8pt](0.5,2){.5}{90}{180}
 \psarc[linewidth=.8pt](1.5,2){.5}{0}{90}
 \psline[linewidth=.8pt](0.5,2.5)(1.5,2.5)

 \psset{origin={0,0}}
 \psline{->}(-.5,1)(.5,1)
 \rput(-7.2,1){$\beta=$}
 \rput(6.5,1){$=\beta\Upsilon$}

 \caption{The action of the map $\Upsilon:\varepsilon\mathscr{I}_n^*\varepsilon\to\mathscr{I}_{n-1}^*$ in the case $n=5$.}\label{Upsilon}
   \end{center}
 \end{figure}

\subsection{Presentations}

Let $X$ be an alphabet (a set whose elements are called \emph{letters}), and
denote by $X^*$ the free monoid on $X$. For $R\subseteq X^*\times X^*$ we denote by
$R^\sharp$ the congruence on $X^*$ generated by~$R$, and we define
$\langle X|R\rangle =X^*/R^\sharp$. We say that a monoid $M$ \emph{has presentation}
$\langle X|R\rangle$ if~${M\cong\langle X|R\rangle}$. Elements of $X$ and $R$ are called
\emph{generators} and \emph{relations} (respectively), and a relation
$(w_1,w_2)\in R$ is conventionally displayed as an equation:~${w_1=w_2}$. We will
often make use of the following universal property of $\langle X|R\rangle$. We say that
a monoid~$S$ \emph{satisfies}~$R$ (or that $R$ \emph{holds in} $S$) via a map
$i_S:X\to S$ if for all~${(w_1,w_2)\in R}$ we have~$w_1i_S^*=w_2i_S^*$ (where
$i_S^*:X^*\to S$ is the natural extension of $i_S$ to $X^*$).
Then~${M=\langle X~|~R\rangle}$ is the monoid, unique up to isomorphism, which is
universal with respect to the property that it satisfies $R$ (via $i_M:x\mapsto
xR^\sharp$); that is, if a monoid $S$ satisfies $R$ via $i_S$, there is a
\emph{unique} homomorphism $\phi:M\to S$ such that $i_M\phi=i_S$:
   \begin{center}
 \vspace{2 cm}
 \psset{xunit= .6 cm,yunit= 1.1 cm}
 \psset{origin={0,0}}
 \uput[r](1,0){$S$}
 \uput[l](-4,2){$X$}
 \uput[r](1,2){$M=\langle X|R\rangle$}
 \psline{->}(-4.2,1.9)(1.1,0.1)
 \psline{->}(-4.2,2)(1.1,2)
 \psline{->}(1.56,1.8)(1.56,.4)
 \rput(-1.2,2.3){\small $i_M$}
 \rput(-1.6,.7){\small $i_S$}
 \rput(2,1.1){\small $\phi$}
   \end{center}
This map $\phi$ is called the \emph{canonical homomorphism}. If
$X$ generates $S$ via $i_S$, then $\phi$ is surjective since
$i_S^*$ is.

\section{Inverse Monoids Generated by Completely Regular Elements}

In this section we present two general results which give
necessary and sufficient conditions for a monoid generated by
completely regular elements to be inverse, with a semilattice of idempotents specified by the generators.

For a monoid $S$, we write $E(S)$ and $G(S)$ for the set of
idempotents and group of units of $S$ (respectively). Suppose now
that $S$ is an inverse monoid (so that $E(S)$ is in fact a
semilattice). The \emph{factorizable part} of $S$ is
$F(S)=E(S)G(S)=G(S)E(S)$, and $S$ is \emph{factorizable} if
$S=F(S)$; in general, $F(S)$ is the largest factorizable inverse
submonoid of~$S$.

Recall that an element $x$ of a monoid $S$ is said to be \emph{completely regular}
if its $\mathscr H$-class $H_x$ is a group. For a completely regular element $x\in S$, we
write $x^{-1}$ for the inverse of $x$ in $H_x$, and $x^0$ for the identity element
of $H_x$. Thus,~${xx^{-1}=x^{-1}x=x^0}$ and, of course, $x^0\in E(S)$. If $X\subseteq
S$, we write $X^0=\left\{  x^{0}~|~x\in X\right\}$.

\begin{proposition}\label{secondprop}
Let $S$ be a monoid, and suppose that $S=\langle X\rangle$ with each $x\in
X$ completely regular. Then $S$ is inverse with $E(S)=\langle X^0\rangle$
if and only if, for all $x,y\in X$, 
\bit %
\item[\emph{(i)}] $x^0y^0=y^0x^0$, and %
\item[\emph{(ii)}] $y^{-1}x^0y\in\langle X^0\rangle$. %
\eit
\end{proposition}

\pf If $S$ is inverse, then (i) holds. Also, for $x,y\in X$, we
have
\[
(y^{-1}x^0y)^2 = y^{-1} x^0y^0x^0y = y^{-1}x^0x^0y^0y =
y^{-1}x^0y,
\]
so that $y^{-1}x^0y\in E(S)$. So if $E(S)=\langle X^0\rangle$, then (ii)
holds.

Conversely, suppose that (i) and (ii) hold. From (i) we see that
$\langle X^0\rangle \subseteq E(S)$ and that~$\langle X^0\rangle$ is a semilattice. Now
let $y\in X$. Next we demonstrate, by induction on $n$, that
\begin{gather} 
\tag{A} 
y^{-1}(X^0)^ny \subseteq \langle X^0\rangle, 
\end{gather}
for all $n\in\mathbb N$. Clearly (A) holds for $n=0$. Suppose next
that (A) holds for some $n\in\mathbb N$, and that
$w\in(X^0)^{n+1}$. So $w=x^0v$ for some $x\in X$ and
$v\in(X^0)^n$. But then by (i) we have
\[
y^{-1}wy = y^{-1}x^0vy = y^{-1}y^0x^0vy = y^{-1}x^0y^0vy =
(y^{-1}x^0y )(y^{-1}vy).
\]
Condition (ii) and an inductive hypothesis that $y^{-1}vy\in\langle
X^0\rangle$ then imply ${y^{-1}wy\in\langle X^0\rangle}$. So (A) holds.

Next we claim that for each $w\in S$ there exists $w'\in S$ such
that
\begin{gather}
\tag{B1} w'\langle X^0\rangle w \subseteq \langle X^0\rangle,\\
\tag{B2} ww',w'w \in \langle X^0\rangle,\\
\tag{B3} ww'w = w ,\, w'ww'=w'.
\end{gather}
We prove the claim by induction on the \emph{length} of $w$ (that is, the minimal
value of $n\in\mathbb N$ for which $w\in X^n$). The case $n=0$ is trivial since
then $w=1$ and we may take $w'=1$. Next suppose that $n\in\mathbb N$ and that the
claim is true for elements of length $n$. Suppose that~$w\in S$ has length $n+1$,
so that $w=xv$ for some $x\in X$ and $v\in S$ of length $n$. Put~$w'=v'x^{-1}$.
Then
\[
w'\langle X^0\rangle w = v'x^{-1}\langle X^0\rangle xv \subseteq v'\langle X^0\rangle v \subseteq
\langle X^0\rangle,
\]
the first inclusion holding by (A) above, and the second by
inductive hypothesis. Thus (B1) holds. Also,
\[
ww' = xvv'x^{-1} \in x\langle X^0\rangle x^{-1} \subseteq \langle X^0\rangle 
\]
and
\[ w'w =
v'x^{-1}xv = v'x^0v \in v'\langle X^0\rangle v \subseteq \langle X^0\rangle
\]
by (A), (B1), and the induction hypothesis, establishing (B2). For (B3), we have
\[
ww'w = xvv'x^{-1}xv = xvv'x^0 v = xx^0vv'v = xv = w,
\]
using (B2), (i), and the inductive hypothesis. Similarly we have
\[
w'ww' = v'x^0vv'x^{-1} = v'vv'x^0x^{-1} = v'x^{-1} =w',
\]
completing the proof of (B3).

Since $S$ is regular, by (B3), the proof will be complete if we can show that
$E(S)\subseteq\langle X^0\rangle$. So suppose that $w\in E(S)$, and choose $w'\in S$ for which
(B1---B3) hold. Then
\[
w' = w'ww' = (w'w)(ww') \in \langle X^0\rangle
\]
by (B2), whence $w'\in E(S)$. But then
\[
w = ww'w = (ww')(w'w) \in\langle X^0\rangle,
\]
again by (B2). This completes the proof. \epf

\begin{proposition}\label{thirdprop}
Suppose that $S$ is a monoid and that $S=\langle G\cup\{z\}\rangle$ where $G=G(S)$
and~$z^3=z$. Then $S$ is inverse with
\[
E(S)=\langle {g^{-1}z^2g~|~g\in G}\rangle \text{~~and~~} F(S)=\langle G\cup\{z^2\} \rangle
\]
if and only if, for all $g\in G$,
\begin{align}
\tag{C1} g^{-1}z^2gz^2 &= z^2g^{-1}z^2g \\
\tag{C2} zg^{-1}z^2gz &\in \langle G\cup\{z^2\}\rangle.
\end{align}
\end{proposition}

\pf First observe that $z$ is completely regular, with $z=z^{-1}$
and $z^0=z^2$. Now put
\[
X=G\cup \{g^{-1}zg~|~g\in G\}.
\]
Then $S=\langle X\rangle$, and each $x\in X$ is completely regular.
Further, if $y=g^{-1}zg$ (with $g\in G$), then $y^{-1}=y$ and
$y^0=y^2=g^{-1}z^2g$. Thus,~${X^0 = \{1\}\cup\{g^{-1}z^2g~|~g\in G\}}$.

Now if $S$ is inverse, then (C1) holds. Also, $zg^{-1}z^2gz\in
E(S)$ for all $g\in G$ so that (C2) holds if $F(S)=\langle
G\cup\{z^2\} \rangle$.

Conversely, suppose now that (C1) and (C2) hold. We wish to verify
Conditions (i) and~(ii) of Proposition \ref{secondprop}, so let
$x,y\in X$. If $x^0=1$ or $y^0=1$, then (i) is immediate, so
suppose~${x^0=g^{-1}z^2g}$ and $y^0=h^{-1}z^2h$ (where $g,h\in
G$). By (C1) we have
\[
x^0y^0 = h^{-1}(hg^{-1}z^2gh^{-1})z^2h =
h^{-1}z^2(hg^{-1}z^2gh^{-1})h = y^0x^0,
\]
and (i) holds.

If $x^0=1$ or $y\in G$, then (ii) is immediate, so suppose
$x^0=g^{-1}z^2g$ and $y=h^{-1}zh$ (where~$g,h\in G$). Then
$y=y^{-1}$ and, by (C2),
\[
y^{-1}x^0y = h^{-1}(zh g^{-1}z^2g h^{-1}z) h \in h^{-1} \langle
G\cup\{z^2\} \rangle h \subseteq\langle G\cup\{z^2\} \rangle.
\]
But by \cite[Lemma 2]{1} and (C1), $\langle G\cup\{z^2\} \rangle$ is a
factorizable inverse submonoid of $S$ with~${E\big(\langle
G\cup\{z^2\} \rangle\big)=\langle X^0\rangle}$. Since
\[
(y^{-1}x^0y)^2 = y^{-1}x^0y^0x^0y = y^{-1}x^0y \in E\big(\langle
G\cup\{z^2\} \rangle\big),
\]
it follows that $y^{-1}x^0y\in\langle X^0\rangle$, so that (ii) holds. So, by Proposition
\ref{secondprop}, $S$ is inverse with 
$E\left( S\right) =\left\langle X^{0}\right\rangle =\left\langle
g^{-1}z^{2}g~|~g\in G\right\rangle $ and, moreover, its factorizable part satisfies
\[
F(S) =  E(S)G \subseteq \langle G\cup X^0\rangle \subseteq \langle G\cup\{z^2\} \rangle \subseteq
F(S).
\]
Hence $F(S)=\langle G\cup\{z^2\}\rangle$, and the proof is complete. \epf

\section{A Presentation of $\mathscr{I}_n^*$}

If $n\leq2$, then $\mathscr{I}_n^*=\mathscr{F}_n$ is equal to its factorizable part. A presentation
of $\mathscr{F}_n$ (for any~$n$) may be found in \cite{1} so, without loss of
generality, we will assume for the remainder of the article that $n\geq3$. We
first fix an alphabet
\[
\mathscr{X}=\mathscr{X}_n=\{x,s_1,\ldots,s_{n-1}\}.
\]
Several notational conventions will prove helpful, and we note
them here. The empty word will be denoted by $1$. A word
$s_i\cdots s_j$ is assumed to be empty if either (i) $i>j$ and the
subscripts are understood to be ascending, or (ii) if $i<j$ and
the subscripts are understood to be descending.

For $1\leq i,j\leq n-1$, we define integers
\[
m_{ij} = \begin{cases} %
1 &\quad\text{if\, $i=j$}\\
3 &\quad\text{if\, $|i-j|=1$}\\
2 &\quad\text{if\, $|i-j|>1$.}
\end{cases}
\]
It will be convenient to use abbreviations for certain words in
the generators which will occur frequently in relations and
proofs. Namely, we write
\[
\sigma = s_2s_3s_1s_2, 
\]
and inductively we define words $l_2,\ldots,l_{n-1}$ and $y_3,\ldots,y_n$ by
\begin{align*}
l_2=xs_2s_1 &\AND l_{i+1}=s_{i+1}l_is_{i+1}s_i &\hspace{-1 cm}\text{for ${2\leq
i\leq
n-2}$,}\\
\intertext{and} %
y_3=x &\AND y_{i+1}=l_iy_is_i &\hspace{-2 cm}\text{for $3\leq i\leq n-1$.}
\end{align*}

Consider now the set $\mathscr{R}=\mathscr{R}_n$ of relations
\begin{align}
\tag{R1} (s_is_j)^{m_{ij}} &= 1 &&\text{for\, $1\leq i\leq
j\leq n-1$}\\
\tag{R2} x^3 &= x \\
\tag{R3} xs_1=s_1x &= x \\
\tag*{} xs_2x =xs_2xs_2 &= s_2xs_2x\\
\tag{R4} &=  xs_2x^2=x^2s_2x\\
\tag{R5} x^2\sigma x^2\sigma = \sigma x^2\sigma x^2 &= xs_2s_3s_2x \\
\tag{R6} y_is_iy_i &= s_iy_is_i &&\text{for\, $3\leq i\leq n-1$} \\
\tag{R7} xs_i &= s_ix &&\text{for\, $4\leq i\leq n-1$.}
\end{align}
Before we proceed, some words of clarification are in order. We
say a relation belongs to~$\mathscr{R}_n$ \emph{vacuously} if it involves a
generator $s_i$ which does not belong to $\mathscr{X}_n$; for example,~(R5)
is vacuously present if $n=3$ because $\mathscr{X}_3$ does not contain the
generator $s_3$. So the reader might like to think of $\mathscr{R}_n$ as
the set of relations (R1---R4) if $n=3$, (R1---R6) if $n=4$, and
(R1---R7) if $n\geq5$. We also note that we will mostly refer only
to the $i=3$ case of relation (R6), which simply says
$xs_3x=s_3xs_3$.

We aim to show that $\mathscr{I}_n^*$ has presentation $\langle \mathscr{X}~|~\mathscr{R} \rangle$, so put
$M=M_n=\langle \mathscr{X}~|~\mathscr{R} \rangle=\mathscr{X}^*/\mathscr{R}^\sharp$. Elements of $M$ are $\mathscr{R}^\sharp$-classes of words over
$\mathscr{X}$. However, in order to avoid cumbersome notation, we will think of elements of
$M$ simply as words over $\mathscr{X}$, identifying two words if they are equivalent under
the relations $\mathscr{R}$. Thus, the reader should be aware of this when reading
statements such as ``\emph{Let $w\in M$}'' and so on.

With our goal in mind, consider the map
\[
\Phi=\Phi_n:\mathscr{X}\to\mathscr{I}_n^*
\]
defined by %

$$\textstyle{\text{$x\Phi = \big( {1,2 \atop 3} \big|
{3 \atop 1,2} \big| {4 \atop 4} \big| {\cdots \atop \cdots} \big|
{n \atop n} \big)$ \ \ and \ \ $s_i\Phi = \big( {1 \atop 1}
\big|{\cdots \atop \cdots} \big|{i-1 \atop i-1} \big|{i \atop i+1}
\big|{i+1 \atop i} \big|{i+2 \atop i+2} \big|{\cdots \atop
\cdots}\big|{n \atop n} \big)$ for $1\leq i\leq n-1$.}}$$ See also
Fig. 4 for illustrations.

\begin{figure}[ht]
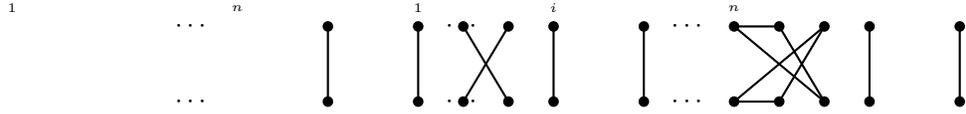

   \begin{center}
 \vspace{1.3 cm}
 \psset{xunit= .6 cm,yunit=.5 cm}
 \psset{origin={1,0}}
 \uvertex0
 \uvertex2
 \uvertex3
 \uvertex4
 \uvertex5
 \uvertex7
 \lvertex0
 \lvertex2
 \lvertex3
 \lvertex4
 \lvertex5
 \lvertex7
 \stline00
 \stline22
 \stline34
 \stline43
 \stline55
 \stline77
 \psset{origin={-8,0}}
 \uvertex0
 \uvertex1
 \uvertex2
 \uvertex3
 \uvertex5
 \lvertex0
 \lvertex1
 \lvertex2
 \lvertex3
 \lvertex5
 \stline02
 \stline12
 \stline20
 \stline21
 \stline33
 \stline55
 \psline(0,0)(1,0)
 \psline(0,2)(1,2)
 \rput(-4,2){$\cdots$}
 \rput(-4,0){$\cdots$}
 \rput(2,2){$\cdots$}
 \rput(2,0){$\cdots$}
 \rput(7,2){$\cdots$}
 \rput(7,0){$\cdots$}
 \rput(-8,2.5){\tiny $1$}
 \rput(-3,2.5){\tiny $n$}
 \rput(1,2.5){\tiny $1$}
 \rput(4,2.5){\tiny $i$}
 \rput(8,2.5){\tiny $n$}
 \caption{The block bijections $x\Phi$ (left) and $s_i\Phi$ (right) in $\mathscr{I}_n^*$.}
 \label{picofgens}
   \end{center}
 \end{figure}

\begin{lemma}\label{relslemma}
The monoid $\mathscr{I}_n^*$ satisfies $\mathscr{R}$ via $\Phi$.
\end{lemma}

\pf This lemma may be proved by considering the relations
one-by-one and diagrammatically verifying that they each hold.
This is straightforward in most cases, but we include a proof for
the more technical relation (R6). First, one may check that $l_i\Phi$ ($2\leq i\leq n-1$) and
$y_j\Phi$ ($3\leq j\leq n$) have graphical representations as
pictured in Fig. 5.

\begin{figure}[ht]
   \begin{center}
 \vspace{1.5 cm}
 \psset{xunit= .6 cm,yunit=.5 cm}
 \psset{origin={-10,0}}
 \uvertex0
 \uvertex1
 \uvertex2
 \uvertex4
 \uvertex5
 \uvertex6
 \uvertex8
 \lvertex0
 \lvertex1
 \lvertex3
 \lvertex4
 \lvertex5
 \lvertex6
 \lvertex8
 \stline00
 \stline10
 \stline21
 \stline43
 \stline54
 \stline55
 \stline66
 \stline88
 \psline(0,2)(1,2)
 \psline(4,0)(5,0)
 \psset{origin={2,0}}
 \uvertex0
 \uvertex1
 \uvertex3
 \uvertex4
 \uvertex5
 \uvertex6
 \uvertex8
 \lvertex0
 \lvertex1
 \lvertex3
 \lvertex4
 \lvertex5
 \lvertex6
 \lvertex8
 \stline05
 \stline45
 \stline50
 \stline54
 \stline66
 \stline88
 \psline(0,0)(1.5,0)
 \psline(0,2)(1.5,2)
 \psline[linestyle=dotted](1.5,0)(2.5,0)
 \psline[linestyle=dotted](1.5,2)(2.5,2)
 \psline(2.5,0)(4,0)
 \psline(2.5,2)(4,2)
 \rput(-7,2){$\cdots$}
 \rput(-8,0){$\cdots$}
 \rput(-3,2){$\cdots$}
 \rput(-3,0){$\cdots$}
 \rput(4.08,1.95){$\cdots$}
 \rput(4.08,-.05){$\cdots$}
 \rput(9,2){$\cdots$}
 \rput(9,0){$\cdots$}
 \rput(-10,2.5){\tiny $1$}
 \rput(-6,2.5){\tiny $i$}
 \rput(-2,2.5){\tiny $n$}
 \rput(2,2.5){\tiny $1$}
 \rput(7,2.5){\tiny $j$}
 \rput(10,2.5){\tiny $n$}
 \caption{The block bijections $l_i\Phi$ (left) and $y_j\Phi$ (right) in $\mathscr{I}_n^*$.}
 \label{picofliyj}
   \end{center}
 \end{figure}

Using this, we demonstrate in Fig. 6 that relation (R6) holds. \epf

\begin{figure}[ht]
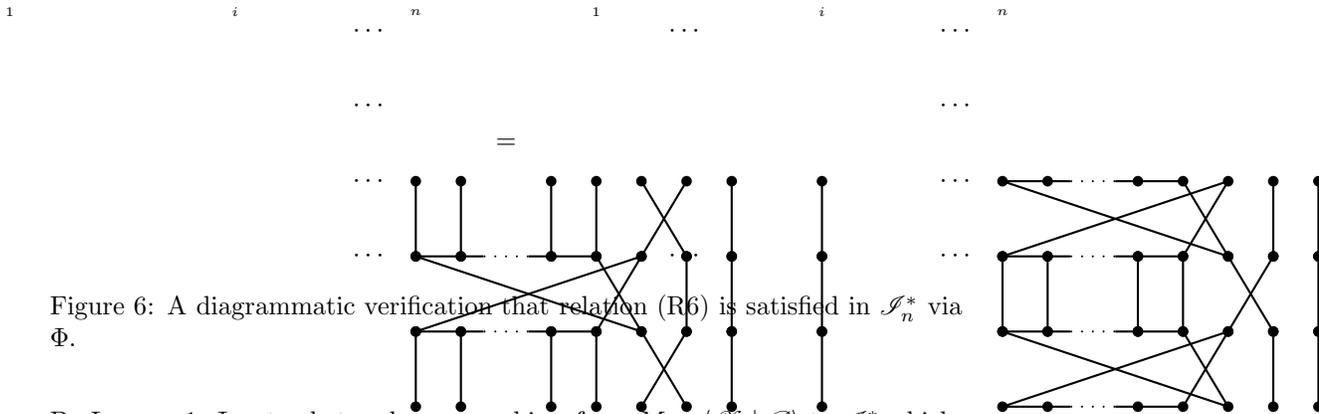

   \begin{center}
 \vspace{3.3 cm}
 \psset{xunit= .6 cm,yunit=.5 cm}
 \psset{origin={-11,0}}
 \uvertex0
 \uvertex1
 \uvertex3
 \uvertex4
 \uvertex5
 \uvertex6
 \uvertex7
 \uvertex9
 \lvertex0
 \lvertex1
 \lvertex3
 \lvertex4
 \lvertex5
 \lvertex6
 \lvertex7
 \lvertex9
 \stline05
 \stline45
 \stline50
 \stline54
 \stline66
 \stline77
 \stline99
 \psline(0,0)(1.5,0)
 \psline(0,2)(1.5,2)
 \psline[linestyle=dotted](1.5,0)(2.5,0)
 \psline[linestyle=dotted](1.5,2)(2.5,2)
 \psline(2.5,0)(4,0)
 \psline(2.5,2)(4,2)
 \psset{origin={-11,4}}
 \uvertex0
 \uvertex1
 \uvertex3
 \uvertex4
 \uvertex5
 \uvertex6
 \uvertex7
 \uvertex9
 \lvertex0
 \lvertex1
 \lvertex3
 \lvertex4
 \lvertex5
 \lvertex6
 \lvertex7
 \lvertex9
 \stline05
 \stline45
 \stline50
 \stline54
 \stline66
 \stline77
 \stline99
 \psline(0,0)(1.5,0)
 \psline(0,2)(1.5,2)
 \psline[linestyle=dotted](1.5,0)(2.5,0)
 \psline[linestyle=dotted](1.5,2)(2.5,2)
 \psline(2.5,0)(4,0)
 \psline(2.5,2)(4,2)
 \psset{origin={2,2}}
 \uvertex0
 \uvertex1
 \uvertex3
 \uvertex4
 \uvertex5
 \uvertex6
 \uvertex7
 \uvertex9
 \lvertex0
 \lvertex1
 \lvertex3
 \lvertex4
 \lvertex5
 \lvertex6
 \lvertex7
 \lvertex9
 \stline05
 \stline45
 \stline50
 \stline54
 \stline66
 \stline77
 \stline99
 \psline(0,0)(1.5,0)
 \psline(0,2)(1.5,2)
 \psline[linestyle=dotted](1.5,0)(2.5,0)
 \psline[linestyle=dotted](1.5,2)(2.5,2)
 \psline(2.5,0)(4,0)
 \psline(2.5,2)(4,2)
 \psset{origin={-11,2}}
 \uvertex0
 \uvertex1
 \uvertex3
 \uvertex4
 \uvertex5
 \uvertex6
 \uvertex7
 \uvertex9
 \lvertex0
 \lvertex1
 \lvertex3
 \lvertex4
 \lvertex5
 \lvertex6
 \lvertex7
 \lvertex9
 \stline00
 \stline11
 \stline33
 \stline44
 \stline56
 \stline65
 \stline77
 \stline99
 \psset{origin={2,0}}
 \uvertex0
 \uvertex1
 \uvertex3
 \uvertex4
 \uvertex5
 \uvertex6
 \uvertex7
 \uvertex9
 \lvertex0
 \lvertex1
 \lvertex3
 \lvertex4
 \lvertex5
 \lvertex6
 \lvertex7
 \lvertex9
 \stline00
 \stline11
 \stline33
 \stline44
 \stline56
 \stline65
 \stline77
 \stline99
 \psset{origin={2,4}}
 \uvertex0
 \uvertex1
 \uvertex3
 \uvertex4
 \uvertex5
 \uvertex6
 \uvertex7
 \uvertex9
 \lvertex0
 \lvertex1
 \lvertex3
 \lvertex4
 \lvertex5
 \lvertex6
 \lvertex7
 \lvertex9
 \stline00
 \stline11
 \stline33
 \stline44
 \stline56
 \stline65
 \stline77
 \stline99
 \rput(-11,6.5){\tiny $1$}
 \rput(-6,6.5){\tiny $i$}
 \rput(-2,6.5){\tiny $n$}
 \rput(2,6.5){\tiny $1$}
 \rput(7,6.5){\tiny $i$}
 \rput(11,6.5){\tiny $n$}
 \rput(-3,0){$\cdots$}
 \rput(-3,2){$\cdots$}
 \rput(-3,4){$\cdots$}
 \rput(-3,6){$\cdots$}
 \rput(10,0){$\cdots$}
 \rput(10,2){$\cdots$}
 \rput(10,4){$\cdots$}
 \rput(10,6){$\cdots$}
 \rput(4,0){$\cdots$}
 \rput(4,6){$\cdots$}
 \rput(0,3){$=$}
 \caption{A diagrammatic verification that relation (R6) is satisfied in $\mathscr{I}_n^*$ via $\Phi$.}
 \label{picofR6}
   \end{center}
 \end{figure}

By Lemma \ref{relslemma}, $\Phi$ extends to a homomorphism from $M=\langle \mathscr{X}~|~\mathscr{R} \rangle$ to
$\mathscr{I}_n^*$ which, without causing confusion, we will also denote by $\Phi=\Phi_n$.
By \cite[Proposition 16]{5}, $\mathscr{I}_n^*$ is generated by~$\mathscr{X}\Phi$, so that $\Phi$
is in fact an epimorphism. Thus, it remains to show that $\Phi$ is injective, and
the remainder of the paper is devoted to this task. The proof we offer is perhaps
unusual in the sense that it uses, in the general case, not a normal form for elements of $M$, but rather structural information about $M$ and an inductive argument. The induction is founded on the case~${n=3}$, for which a normal form is given in the next proposition.

\begin{proposition}\label{n=3case}
The map $\Phi_3$ is injective.
\end{proposition}

\pf Consider the following list of 25 words in $M_3$: %
\begin{itemize}
\item the 6 units $\{1,s_1,s_2,s_1s_2,s_2s_1,s_1s_2s_1\}$, %
\item the 18 products in
$\{1,s_2,s_1s_2\}\{x,x^2\}\{1,s_2,s_2s_1\}$, and %
\item the zero element $xs_2x$.
\end{itemize} %
This list contains the generators, and is easily checked to be closed under
multiplication on the right by the generators. Thus, $|M_3|\leq 25$. But $\Phi_3$ is a
surjective map from $M_3$ onto $\mathscr{I}_3^*$, which has cardinality $25$. It follows
that $|M_3|=25$, and that $\Phi_3$ is injective. \epf

From this point forward, we assume that $n\geq4$.
The inductive step in our argument relies on Proposition \ref{firstprop} below,
which provides a sufficient condition for a homomorphism of inverse monoids to be
injective. Let $S$ be an inverse monoid and, for $s,t\in S$, write~$s^t=t^{-1}st$.
We say that a non-identity idempotent $e\in E(S)$ has property (P) if, for all
non-identity idempotents $f\in E(S)$, there exists $g\in G(S)$ such that~$f^g\in
eSe$.

\begin{proposition}\label{firstprop}
Let $S$ be an inverse monoid with $E=E(S)$ and $G=G(S)$. Suppose that~$1\not=e\in
E$ has property (P). Let $\phi:S\to T$ be a homomorphism of inverse monoids for
which $\phi|_E$, $\phi|_G$, and $\phi|_{eSe}$ are injective. Then $\phi$ is
injective.
\end{proposition}

\pf By the kernel-and-trace description of congruences on $S$ \cite[Section
5.1]{4}, and the injectivity of $\phi|_G$, it is enough to show that
$x\phi=f\phi$ (with $x\in S$ and $1\not=f\in E$) implies $x=f$, so suppose that
$x\phi=f\phi$. Choose $g\in G$ such that $f^g\in eSe$. Now $
(xx^{-1})\phi=f\phi=(x^{-1}x)\phi, $ so that $f=xx^{-1}=x^{-1}x$, since $\phi|_E$
is injective. Thus $f \mathscr{H} x$ and it follows that $f^g \mathscr{H} x^g$, so that $x^g\in
eSe$.  Now $x\phi=f\phi$ also implies $x^g\phi=f^g\phi$ and so, by the injectivity
of $\phi|_{eSe}$, we have $x^g=f^g$, whence $x=f$. \epf

It is our aim to apply Proposition \ref{firstprop} to the map $\Phi:M\to\mathscr{I}_n^*$.
In order to do this, we first use Proposition \ref{thirdprop} to show (in Section
\ref{sect:structure}) that~$M$ is inverse, and we also deduce information about
its factorizable part~$F(M)$, including the fact that $\Phi|_{F(M)}$ is injective;
this then implies that both~$\Phi|_{E(M)}$ and~$\Phi|_{G(M)}$ are injective too.
Finally, in Section \ref{sect:local}, we locate a non-identity idempotent $e\in M$
which has property (P). We then show that the injectivity of $\Phi|_{eMe}$ is
equivalent to the injectivity of $\Phi_{n-1}$ which we assume, inductively. We
first pause to make some observations concerning the factorizable part $\mathscr{F}_n$ of
$\mathscr{I}_n^*$.

\subsection{The factorizable part of $\mathscr{I}_n^*$}\label{sect:factIn*}

Define an alphabet $\mathscr{X}_F=\{t,s_1,\ldots,s_{n-1}\}$, and consider
the set $\mathscr{R}_F$ of relations
\begin{align}
\tag{F1} (s_is_j)^{m_{ij}} &= 1 &&\text{for\, $1\leq i\leq
j\leq n-1$}\\
\tag{F2} t^2 &= t \\
\tag{F3} ts_1=s_1t &=t \\
\tag{F4} ts_i &= s_it &&\text{for\, $3\leq i\leq n-1$}\\
\tag{F5} ts_2ts_2 &= s_2ts_2t\\
\tag{F6} t\sigma t\sigma &= \sigma t\sigma t.
\end{align}
(Recall that $\sigma$ denotes the word $s_2s_3s_1s_2$.) The following
result was proved in \cite{1}.

\begin{thm}\label{Fnpres}
The monoid $\mathscr{F}_n=F(\mathscr{I}_n^*)$ has presentation 
$\langle \mathscr{X}_F~|~\mathscr{R}_F \rangle$ via
\begin{equation*} 
s_i\mapsto s_i\Phi,\,\, t \mapsto (1,2\,|\,3\,|\cdots|\,n). \hfill \Box
\end{equation*} 
\end{thm}

\begin{lemma}\label{RFholdsinM}
The relations $\mathscr{R}_F$ hold in $M$ via the map $\Theta:t\mapsto
x^2,\,s_i\mapsto s_i$.
\end{lemma}

\pf Relations (F1---F3) are immediate from (R1---R3); (F5) follows  from several
applications of~(R4); and (F6) forms part of (R5). The $i\geq4$ case of (F4)
follows from (R7), and the $i=3$ case follows from~(R1) and (R6), since
\begin{equation*}
x^2s_3 = xs_3s_3xs_3 = xs_3 xs_3x = s_3xs_3 s_3x = s_3x^2.
\end{equation*}
\epf

It follows that $\Theta\circ\Phi$ extends to a homomorphism of
$\langle \mathscr{X}_F~|~\mathscr{R}_F \rangle$ to $\mathscr{F}_n$, which is an isomorphism by Theorem \ref{Fnpres}. We conclude that 
$\Phi|_{\langle x^2,s_1,\ldots,s_{n-1}\rangle} = \Phi|_{\text{im} (\Theta)}$ 
is injective (and therefore an isomorphism).

\subsection{The structure of $M$}\label{sect:structure}

It is easy to see that the group of units $G(M)$ is the subgroup generated by
$\{s_1,\ldots,s_{n-1}\}$. The reason for this is that relations (R2---R7) contain
at least one occurrence of $x$ on both sides. Now (R1) forms the set of defining
relations in Moore's famous presentation~\cite{6} of the symmetric group
$\mathcal{S}_n$. Thus we may identify $G(M)$ with~$\mathcal{S}_n$ in the obvious way. Part (i) of
the following well-known result (Lemma 8) gives a normal form for the elements of $\mathcal{S}_n$
(and is probably due to Burnside; a proof is also sketched in \cite{1}). The
second part follows immediately from the first, and is expressed in terms of a
convenient contracted notation which is defined as follows. Let $1\leq i\leq n-1$,
and $0\leq k\leq n-1$. We write
\[
s_i^k = \begin{cases} %
s_i &\quad\text{if\, $i\leq k$}\\
1 &\quad\text{otherwise.}
\end{cases}
\]
The reader might like to think of this as abbreviating $s_i^{k\geq
i}$, where $k\geq i$ is a boolean value, equal to $1$ if
$k\geq i$ holds and $0$ otherwise.

\begin{lemma}\label{Snnorm}
Let $g\in G(M)=\langle s_1,\ldots,s_{n-1}\rangle$. Then %
\begin{description} %
\item[\emph{(i)}] $g=(s_{i_1}\cdots s_{j_1})\cdots(s_{i_k}\cdots
s_{j_k})$ for some $k\geq0$ and some $i_1\leq j_1,\ldots,i_k\leq
j_k$ with $1\leq i_k<\cdots<i_1\leq n-1$, and %
\item[\emph{(ii)}] $g=hs_2^ks_3^ks_4^k(s_5\cdots s_k)s_1^\ell
s_2^\ell s_3^\ell(s_4\cdots s_\ell) = hs_2^ks_3^ks_4^ks_1^\ell
s_2^\ell s_3^\ell(s_5\cdots s_k)(s_4\cdots s_\ell)$ for some
$h\in\langle s_3,\ldots,s_{n-1}\rangle$, $k\geq 1$ and $\ell\geq 0$. \hfill $\Box$ %
\end{description} 
\end{lemma} 

We are now ready to prove the main result of this section.

\begin{proposition}\label{Misinverse}
The monoid $M=\langle \mathscr{X}~|~\mathscr{R} \rangle$ is inverse, and we have
\[
E(M) = \langle g^{-1}x^2g~|~g\in G(M)\rangle  \text{~~~and~~~}  
F(M)=\langle x^2, s_1,\ldots,s_{n-1}\rangle.
\]
\end{proposition}

\pf Put $G=G(M)=\langle s_1,\ldots,s_{n-1}\rangle$. So $M=\langle G\cup\{x\}\rangle$ and $x=x^3$.
We will now verify conditions (C1) and (C2) of Proposition \ref{thirdprop}.
By Lemma \ref{RFholdsinM}, $\langle G\cup\{x^2\}\rangle$ is a homomorphic (in fact
isomorphic) image of $\mathscr{F}_n$, so $g^{-1}x^2g$ commutes with $x^2$ for all $g\in G$,
and condition (C1) is verified.

To prove (C2), let $g\in G$. By Lemma \ref{Snnorm}, we have $$g=
hs_2^ks_3^ks_4^ks_1^\ell s_2^\ell s_3^\ell(s_5\cdots
s_k)(s_4\cdots s_\ell)$$ for some $h\in\langle s_3,\ldots,s_{n-1}\rangle$,
$k\geq 1$ and $\ell\geq 0$. Now
\begin{align*}
&xg^{-1}x^2gx \\&= x(s_\ell\cdots s_4)(s_k\cdots s_5) s_3^\ell
s_2^\ell s_1^\ell s_4^ks_3^ks_2^k
(h^{-1}x^2h)s_2^ks_3^ks_4^ks_1^\ell s_2^\ell s_3^\ell(s_5\cdots
s_k)(s_4\cdots s_\ell)x \\
&= (s_\ell\cdots s_4)(s_k\cdots s_5) xs_3^\ell s_2^\ell s_1^\ell
s_4^ks_3^ks_2^k x^2 s_2^ks_3^ks_4^ks_1^\ell s_2^\ell s_3^\ell
x(s_5\cdots s_k)(s_4\cdots s_\ell),
\end{align*}
by (R7), (F4), and (R1). Thus it suffices to show that $x(x^2)^\pi x\in\langle
G\cup\{x^2\}\rangle$, where we have written $\pi=s_2^ks_3^ks_4^ks_1^\ell s_2^\ell
s_3^\ell$. Altogether there are 16 cases to consider for all pairs~$(k,\ell)$
with~$k=1,2,3,\geq4$ and $\ell=0,1,2,\geq3$. Table 1 below contains an equivalent
form of~$x(x^2)^\pi x$ as a word over $\{x^2,s_1,\ldots,s_{n-1}\}$ for each
$(k,\ell)$, as well as a list of the relations used in deriving the expression. We
performed the calculations in the order determined by going along the first row
from left to right, then the second, third, and fourth rows. Thus, as for example
in the case $(k,\ell)=(2,1)$, we have used expressions from previously considered
cases.
\begin{table}[ht]
{\footnotesize
\begin{center}

\begin{tabular}{|c||c|c|c|c|c}%
\hline %
  &  $\ell=0$  &  $\ell=1$  &  $\ell=2$  &  $\ell\geq 3$\\
\hline %
\hline %
              &  $x^2$  &  $x^2$   &  $x^2s_2x^2$  &  $x^2\sigma x^2\sigma$  \\
$k=1$         &   (R2)  &  (R2,3)  &  (R3,4)       &  (R1,2,3,5)       \\
              &         &          &               &  and (F4)         \\
\hline %
              &  $x^2s_2x^2$       &  $x^2s_2x^2$       &  $x^2s_2x^2$       &  $x^2\sigma x^2\sigma$       \\
$k=2$         &  (R4)              &  (R3,4)            &  (R1,3,4)          &  (R1,3) and            \\
              &                    &                    &                    &  $(k,\ell)=(1,\geq3)$  \\
\hline %
              &  $x^2\sigma x^2\sigma$       &  $x^2\sigma x^2\sigma$       &  $x^2s_2s_3s_2x^2$  &  $x^2s_2s_3s_2x^2$  \\
$k=3$         &  $(k,\ell)=(1,\geq3)$  &  (R3) and              &  (R1,2,5)           &  (R1,3) and         \\
              &                        &  $(k,\ell)=(1,\geq3)$  &                     &  $(k,\ell)=(3,2)$   \\
\hline %
              &  $s_4x^2\sigma x^2\sigma s_4$  &  $s_4x^2\sigma x^2\sigma s_4$  &  $s_4x^2s_2s_3s_2x^2s_4$  &  $s_3s_4x^2\sigma x^2\sigma s_4s_3$  \\
$k\geq4$         &  (R7) and   &  (R3) and   &  (R1,7) and   &  (R1,2,5,6,7)   \\
              &  $(k,\ell)=(1,\geq3)$  &  $(k,\ell)=(\geq4,0)$  &  $(k,\ell)=(3,2)$  &  and (F4)  \\
\hline %
\end{tabular}

\end{center}
\vspace{-5mm}  }
\caption{Expressions for $x(x^2)^\pi x$ and the
relations used. See text for further explanation.}
\label{table1}
\end{table}
\newline In order that readers need not perform all the calculations
themselves, we provide a small number of sample derivations. The
first case we consider is that in which $(k,\ell)=(1,2)$. In this
case we have $\pi=s_1s_2$ and, by (R3) and several applications of
(R4), we calculate
\[
x(x^2)^\pi x = xs_2s_1x^2s_1s_2x = xs_2x^2s_2x = x^2s_2x^2.
\]
Next suppose $(k,\ell)=(1,\geq3)$. (Here we mean that $k=1$ and
$\ell\geq3$.) Then $\pi=s_1s_2s_3$ and
\begin{align*}
x(x^2)^\pi x &= xs_3s_2s_1x^2s_1s_2s_3x \\
&= xs_3s_2x^2s_2s_3x &&\text{by (R3)}\\
&= xs_3s_2s_3x^2s_3s_2s_3x &&\text{by (R1) and (F4)}\\
&= (x^2\sigma x^2\sigma)(\sigma x^2\sigma x^2) &&\text{by (R1) and (R5)}\\
&= x^2\sigma x^2\sigma x^2 &&\text{by (R1) and (R2)}\\
&= x^2\sigma x^2\sigma &&\text{by (R5) and (R2).}
\end{align*}
If $(k,\ell)=(3,2)$, then $\pi=\pi^{-1}=\sigma$ by (R1) and, by (R2) and (R5), we
have
\[
x(x^2)^\pi x = x\sigma x^2\sigma x = x(\sigma x^2 \sigma x^2) x =
x(xs_2s_3s_2x)x.
\]
If $(k,\ell)=(\geq4,2)$ then $\pi=s_2s_3s_4s_1s_2=\sigma s_4$ by
(R1), and so, using (R7) and the $(k,\ell)=(3,2)$ case, we have
\[
x(x^2)^\pi x = xs_4\sigma x^2\sigma s_4x = s_4x\sigma x^2\sigma xs_4 =
s_4(x^2s_2s_3s_2x^2)s_4.
\]
Finally, we consider the $(k,\ell)=(\geq4,\geq3)$ case. Here we
have $\pi=s_2s_3s_4s_1s_2s_3=\sigma s_4s_3$ by (R1), and so
\begin{align*}
x(x^2)^\pi x &= xs_3s_4\sigma x^2\sigma s_4s_3x\\
&= xx^2s_3s_4\sigma x^2\sigma s_4s_3x &&\text{by (R2)}\\
&= xs_3s_4x^2\sigma x^2\sigma s_4s_3x &&\text{by (F4)}\\
&= xs_3s_4xs_2s_3s_2x s_4s_3x &&\text{by (R5)}\\
&= xs_3xs_4s_2s_3s_2s_4x s_3x &&\text{by (R7)}\\
&= s_3xs_3s_4s_3s_2s_3s_4 s_3xs_3 &&\text{by (R1) and (R6)}\\
&= s_3xs_4s_3s_4s_2s_4s_3 s_4xs_3 &&\text{by (R1)}\\
&= s_3s_4xs_3s_2s_3 xs_4s_3 &&\text{by (R1) and (R7)}\\
&= s_3s_4(x^2\sigma x^2\sigma)s_4s_3 &&\text{by (R1) and (R5).}
\end{align*}
After checking the other cases, the proof is complete. \epf

After the proof of Lemma \ref{RFholdsinM}, we observed that $\Phi|_{\langle
x^2,s_1,\ldots,s_{n-1}\rangle}$ is injective. By Proposition \ref{Misinverse}, we
conclude that $\Phi|_{F(M)}$ is injective. In particular, both $\Phi|_{E(M)}$
and~$\Phi|_{G(M)}$ are injective.

\subsection{A local submonoid}\label{sect:local}

Now put $e=x^2\in M$. So clearly $e$ is a non-identity idempotent
of $M$. Our goal in this section is to show that $e$ has property
(P), and that $eMe$ is a homomorphic image of $M_{n-1}$.

\begin{lemma}\label{propP}
The non-identity idempotent $e=x^2\in M$ has property (P).
\end{lemma}

\pf Let $1\not=f\in E(M)$. By Proposition \ref{Misinverse} we have
$f=e^{g_1}e^{g_2}\cdots e^{g_k}$ for some~$k\geq1$ and $g_1,g_2,\ldots,g_k\in
G(M)$. But then 
$$f^{g_1^{-1}}= e\,e^{g_2g_1^{-1}}\cdots e^{g_kg_1^{-1}}e\in eMe.$$
\epf

We now define words
\[
X=s_3x\sigma xs_3 ,\, S_1=x ,\quad\text{and}\quad S_j = es_{j+1}  \quad\text{for\,
$j=2,\ldots,n-2$,}
\]
and write $\mathscr{Y}=\mathscr{Y}_{n-1}=\{X,S_1,\ldots,S_{n-2}\}$. We note that $e$ is a left and
right identity for the elements of $\mathscr{Y}$ so that $\mathscr{Y}\subseteq eMe$, and that $X=y_4$ (by
definition).

\begin{proposition}\label{eMegens}
The submonoid $eMe$ is generated (as a monoid with identity $e$) by $\mathscr{Y}$.
\end{proposition}

\pf We take $w\in M$ with the intention of showing that $u=ewe\in eMe$ belongs
to~$\langle\mathscr{Y}\rangle$. We do this by induction on the (minimum) number $d=d(w)$ of
occurrences of~$x^\delta$~($\delta\in\{1,2\}$) in $w$. Suppose first that $d=0$, so that
$u=ege$ where $g\in G(M)$. By Lemma \ref{Snnorm} we have
\[
g=h(s_2^js_3^js_1^is_2^i)(s_4\cdots s_j)(s_3\cdots s_i)
\]
for some $h\in\langle s_3,\ldots,s_{n-1}\rangle$ and $j\geq1$, $i\geq0$. Put
$h'=(s_4\cdots s_j)(s_3\cdots s_i)\in\langle s_3,\ldots,s_{n-1}\rangle$. Now by (F2) and
(F4) we have
\[
u = ege = eh \cdot e(s_2^js_3^js_1^is_2^i)e \cdot eh'.
\]
By (F2) and (F4) again, we see that $eh,eh'\in\langle S_2,\ldots,S_{n-1}\rangle$, so it is
sufficient to show that the word~${e\pi e}$ belongs to $\langle\mathscr{Y}\rangle$, where we have
written $\pi=s_2^js_3^js_1^is_2^i$.
Table 2 below contains an equivalent form of $e\pi e$ as a word
over $\mathscr{Y}$ for each $(i,j)$, as well as a list of the relations
used in deriving the expression.
\begin{table}[ht!]
{\footnotesize
\begin{center}
\begin{tabular}{|c||c|c|c|c}%
\hline %
  &  $j=1$  &  $j=2$  &  $j\geq3$  \\
\hline %
\hline %
    &  $e$  &  $X^2$  &  $X^2S_2$    \\
$i=0$         &  (R2)  &  (R1,2,5)   &  (R2), (F4), and      \\
         &    &  and (F4)  &   $(i,j)=(0,2)$   \\
\hline %
    &  $e$  &  $X^2$  &  $X^2S_2$    \\
$i=1$         &  (R2,3)  &  (R3) and   &  (R3) and      \\
         &    &  $(i,j)=(0,2)$  &  $(i,j)=(0,\geq3)$    \\
\hline %
    &  $X^2$  &  $X^2$  &  $S_1S_2XS_2S_1$    \\
$i\geq2$         &  (R3) and  &  (R1,3) and   &  (R1,2)     \\
         &  $(i,j)=(0,2)$  &  $(i,j)=(0,2)$  &   and (F4)   \\
\hline %
\end{tabular}
\end{center}
\vspace{-5mm} \caption{Expressions for $e\pi e$ and the relations
used. See text for further explanation.} }
\label{table2}
\end{table}

Most of these derivations are rather straightforward,
but we include two example calculations. For the $(i,j)=(0,2)$
case, note that
\begin{align*}
X^2 = s_3x\sigma xs_3s_3x\sigma x s_3
&= s_3x(x^2\sigma x^2\sigma) x s_3 &&\text{by (R1) and (R2)}\\
&= s_3x(xs_2s_3s_2x) x s_3 &&\text{by (R5)}\\
&= s_3x^2s_3s_2s_3x^2 s_3 &&\text{by (R1)}\\
&= x^2s_2x^2  &&\text{by (F4) and (R1)}\\
&= e\pi e.\\
\intertext{For the $(i,j)=(\geq2,\geq3)$ case, we have $\pi=\sigma$ and}
e\pi e = x^2\sigma x^2
&= xs_3s_3x\sigma xs_3s_3x &&\text{by (R1)}\\
&= x(x^2s_3)s_3x\sigma xs_3(s_3x^2)x &&\text{by (R2)}\\
&= x(x^2s_3)s_3x\sigma xs_3(x^2s_3)x &&\text{by (F4)}\\
&= S_1S_2XS_2S_1.
\end{align*}
This establishes the $d=0$ case. Now suppose $d\geq1$, so that $w=vx^\delta g$ for
some $\delta\in\{1,2\}$, $v\in M$ with $d(v)=d(w)-1$, and $g\in G(M)$. Then by (R2),
\[
u = ewe = (eve)x^\delta(ege).
\]
Now $x^\delta$ belongs to $\langle\mathscr{Y}\rangle$ since $x^\delta$ is equal to $S_1$
(if $\delta=1$) or $e$ (if $\delta=2$). By an induction hypothesis we
have $eve\in\langle\mathscr{Y}\rangle$ and, by the $d=0$ case considered above, we
also have~${ege\in\langle\mathscr{Y}\rangle}$. \epf

The next step in our argument is to prove (in Proposition
\ref{eMerels} below) that the elements of~$\mathscr{Y}_{n-1}$ satisfy the
relations $\mathscr{R}_{n-1}$ via the obviously defined map. Before we do
this, however, it will be convenient to prove the following basic
lemma. If $w\in M$, we write $\text{rev}(w)$ for the word obtained by
writing the letters of $w$ in reverse order. We say that $w$ is
\emph{symmetric} if $w=\text{rev}(w)$.

\begin{lemma}\label{intlem}
If $w\in M$ is symmetric, then $w=w^3$ and $w^2\in E(M)$.
\end{lemma}

\pf Now $z=z^{-1}$ for all $x\in\mathscr{X}$ and it follows that $w^{-1}=\text{rev}(w)$ for all
$w\in M$. So, if~$w$ is symmetric, then $w=w^{-1}$, in which case
$w=ww^{-1}w=w^3$. \epf

\begin{proposition}\label{eMerels} The elements of $\mathscr{Y}_{n-1}$ satisfy
the relations $\mathscr{R}_{n-1}$ via the map
\[
\Psi: \mathscr{X}_{n-1}\to eMe:x\mapsto X,\,s_i\mapsto S_i.
\]
\end{proposition}

\pf We consider the relations from $\mathscr{R}_{n-1}$ one at a time. In order to avoid
confusion, we will refer to the relations from $\mathscr{R}_{n-1}$ as (R1)$'$, (R2)$'$,
etc. We also extend the use of upper case symbols for the element
$\Sigma=S_2S_3S_1S_2$ as well as the words $L_i$ (for~$i=2,\ldots,n-2$) and $Y_j$
(for $j=3,\ldots,n-1$). It will also prove convenient to refer to the idempotents
\[
e_i = e^{(s_2\cdots s_i)(s_1\cdots s_{i-1})}\in E(M),
\]
defined for each $1\leq i\leq n$. Note that $e_1=e$, and that
$e_i\Phi\in\mathscr{I}_n^*$ is the idempotent with domain
$(1\,|\cdots|\,i-1\,|\,i,i+1\,|\,i+2\,|\cdots|\,n)$.

We first consider relation (R1)$'$. We must show that $(S_iS_j)^{m_{ij}}=e$ for
all $1\leq i\leq j\leq n-2$. Suppose first that $i=j$. Now $S_1^2=e$ by definition
and if $2\leq i\leq n-2$ then, by (R1), (R2), (R7), and (F4), we have
$S_i^2=es_{i+1}es_{i+1}=e^2s_{i+1}^2=e$. Next, if $2\leq j\leq n-2$, then
\begin{align*}
(S_1S_j)^{m_{1j}} = (xe&s_{j+1})^{m_{1j}} = (xs_{j+1})^{m_{1j}} \\
&= \begin{cases} %
xs_3xs_3xs_3 = s_3xs_3s_3xs_3=s_3x^2s_3=x^2=e &\quad\text{if\,
$j=2$}\\
xs_{j+1}xs_{j+1}=x^2s_{j+1}^2=x^2=e &\quad\text{if\, $j\geq3$,}
\end{cases}
\end{align*}
by (R1), (R2), (R6), (R7), and (F4). Finally, if $2\leq i\leq j\leq n-2$, then
\[
(S_iS_j)^{m_{ij}} = (es_{i+1}es_{j+1})^{m_{i+1,j+1}} =
e(s_{i+1}s_{j+1})^{m_{i+1,j+1}} = e,
\]
by (R1), (R2), and (F4). This completes the proof for (R1)$'$.

For (R2)$'$, we have
\begin{align*}
X^3 &= (s_3x\sigma xs_3)(s_3x\sigma xs_3)(s_3x\sigma xs_3)\\
&= s_3x\sigma x^2\sigma x^2\sigma xs_3 &&\text{by (R1)}\\
&= s_3xx^2\sigma x^2\sigma \sigma xs_3 &&\text{by (R5)}\\
&= s_3x\sigma  xs_3 &&\text{by (R1) and (R2)}\\
&= X.\\
\intertext{For (R3)$'$, first note that}
XS_1 &=s_3x\sigma xs_3x\\
&=s_3x\sigma s_3xs_3 &&\text{by (R6)}\\
&=s_3xs_1\sigma xs_3 &&\text{by (R1)}\\
&=s_3x\sigma xs_3 &&\text{by (R3)}\\
&=X.
\end{align*}
(Here, and later, we use the fact that $\sigma s_3=s_1\sigma$, and so also $\sigma s_1=s_3\sigma$. These are easily checked using (R1) or by drawing pictures.) The relation $S_1X=X$ is proved by a symmetrical argument.

To prove (R4)$'$, we need to show that $XS_2X$ is a (left and right) zero for $X$
and $S_2$. Since~$X,S_2\in\langle x,s_1,s_2,s_3\rangle$, it suffices to show that $XS_2X$
is a zero for each of $x,s_1,s_2,s_3$. In order to contract the proof, it will be
convenient to use the following ``arrow notation''. If~$a$ and $b$ are elements of
a semigroup, we write $a\larr b$ and $a\rarr b$ to denote the relations~$ab=a$
($a$ is a left zero for $b$) and $ba=a$ ($a$ is a right zero for $b$),
respectively. The arrows may be superimposed, so that $a\lrarr b$ indicates the
presence of both relations. We first calculate
\begin{align*}
XS_2X &= (s_3x\sigma xs_3)es_3(s_3x\sigma xs_3)\\
&= s_3x\sigma xs_3x\sigma xs_3 &&\text{by (R1) and (R2)}\\
&= s_3x\sigma s_3xs_3\sigma xs_3 &&\text{by (R6)}\\
&= s_3xs_1\sigma x\sigma s_1xs_3 &&\text{by (R1)}\\
&= s_3x\sigma x\sigma xs_3 &&\text{by (R3).}
\end{align*}
Put $w=s_3x\sigma x\sigma xs_3$. We see immediately that $w\larr x$ since
\[
wx=s_3x\sigma x\sigma xs_3x = s_3x\sigma x\sigma s_3xs_3 = s_3x\sigma xs_1\sigma xs_3 = s_3x\sigma x\sigma
xs_3 = w,
\]
by (R1), (R3), and (R6), and a symmetrical argument shows that $w\rarr x$. Next,
note that~$w$ is symmetric so that $w=w^3$ and $w\lrarr w^2$, by Lemma
\ref{intlem}. Since $\hspace{.2ex}\lrarr\hspace{.2ex}$ is transitive, the proof
of~(R4)$'$ will be complete if we can show that $w^2\lrarr s_1,s_2,s_3$. Now by
Lemma \ref{intlem} again we have $w^2\in E(M)$ and, since
$w^2\Phi=(e_1e_2e_3)\Phi$ as may easily be checked diagrammatically, we have
$w^2=e_1e_2e_3$ by the injectivity of $\Phi|_{E(M)}$. But $(e_1e_2e_3)\Phi\lrarr
s_1\Phi,s_2\Phi,s_3\Phi$ in $\mathscr{F}_n$ and so, by the injectivity of $\Phi|_{F(M)}$,
it follows that $w^2=e_1e_2e_3\lrarr s_1,s_2,s_3$.

Relations (R5---R7)$'$ all hold vacuously if $n=4$, so for the remainder of the
proof we assume that $n\geq5$.

Next we consider (R5)$'$. Now, by (R2) and (F4), we see that
\[
\Sigma=S_2S_3S_1S_2 = (es_3)(es_4)x(es_3) = s_3s_4xs_3.
\]
In particular, $\Sigma$ is symmetric, by (R7), and $\Sigma^{-1}=\Sigma$. Also, $X=s_3x\sigma
xs_3$ is symmetric and so~$X^2$ is idempotent by Lemma \ref{intlem}. It follows
that $X^2\Sigma X^2\Sigma$ and $\Sigma X^2\Sigma X^2$ are both idempotent. Since $(X^2\Sigma
X^2\Sigma)\Phi=(\Sigma X^2\Sigma X^2)\Phi$, as may easily be checked diagrammatically, we
conclude that $X^2\Sigma X^2\Sigma=\Sigma X^2\Sigma X^2$, by the injectivity of
$\Phi|_{E(M)}$. It remains to check that $XS_2S_3S_2X=\Sigma X^2\Sigma X^2$. Since
$(XS_2S_3S_2X)\Phi=(\Sigma X^2\Sigma X^2)\Phi$, it suffices to show that~$XS_2S_3S_2X\in E(M)$. By (R1), (R2), and (F4), 
\[
XS_2S_3S_2X = (s_3x\sigma xs_3)es_3es_4es_3(s_3x\sigma xs_3) = s_3(x\sigma
x)s_4(x\sigma x)s_3,
\]
so it is enough to show that $v=(x\sigma x)s_4(x\sigma x)$ is
idempotent. We see that
\begin{align*}
v^2 &= x\sigma xs_4x\sigma x^2\sigma xs_4x\sigma x\\
&= x\sigma xs_4x(x^2\sigma x^2\sigma) xs_4x\sigma x &&\text{by (R2)}\\
&= x\sigma xs_4x(xs_2s_3s_2x) xs_4x\sigma x &&\text{by (R5)}\\
&= x\sigma xs_4s_2s_3s_2 s_4x\sigma x &&\text{by (R2) and (R7).}
\end{align*}
Put $u=x\sigma x$. Since $u$ is symmetric, Lemma \ref{intlem} says that $u^2\in
E(M)$. One verifies easily that $(u^2s_4s_2s_3s_2s_4u^2)\Phi=(u^2s_4u^2)\Phi$ in
$\mathscr{F}_n$ and it follows, by the injectivity of~$\Phi|_{F(M)}$, that
$u^2s_4s_2s_3s_2s_4u^2=u^2s_4u^2$. By Lemma \ref{intlem} again, we also have
$u=u^3$ so that
\[
v^2 = us_4s_2s_3s_2 s_4u = u u^2s_4s_2s_3s_2 s_4u^2u = u
u^2s_4u^2u = u s_4u = v,
\]
and (R5)$'$ holds.

Now we consider (R6)$'$, which says $Y_iS_iY_i=S_iY_iS_i$ for $i\geq3$. So we must
calculate the words $Y_i$ which, in turn, are defined in terms of the words $L_i$.
Now ${L_2=XS_2S_1}$ and ${L_{i+1}=S_{i+1}L_iS_{i+1}S_i}$ for $i\geq2$. A
straightforward induction shows that ${L_i=l_{i+1}e}$ for all~${i\geq2}$. This,
together with the definition of the words $Y_i$ (as $Y_3=X$,
and~${Y_{i+1}=L_iY_iS_i}$ for $i\geq3$) and a simple induction, shows that
$Y_i=y_{i+1}$ for all $i\geq3$. But then for~${3\leq i\leq n-2}$, we have
\[
Y_iS_iY_i = y_{i+1}es_{i+1}y_{i+1} = y_{i+1}s_{i+1}y_{i+1} =
s_{i+1}y_{i+1}s_{i+1} = s_{i+1}ey_{i+1}es_{i+1} = S_iY_iS_i.
\]
Here we have used (R6) and (F4), and the fact, verifiable by a simple induction,
that $y_je=ey_j=y_j$ for all $j$.

Relation (R7)$'$ holds vacuously when $n=5$ so, to complete the proof, suppose
$n\geq6$ and~$i\geq4$. Now $Xe=eX=X$ as we have already observed, and
$Xs_{i+1}=s_{i+1}X$ by (R1) and (R7). Thus $X$ commutes with $es_{i+1}=S_i$. This
completes the proof. \epf

\subsection{Conclusion} 
We are now ready to tie together all the loose ends.

\begin{thm}
The dual symmetric inverse monoid $\mathscr{I}_n^*$ has presentation $\langle \mathscr{X}~|~\mathscr{R} \rangle$
via $\Phi$.
\end{thm}

\pf All that remains is to show that $\Phi=\Phi_n$ is injective.
In Proposition \ref{n=3case} we saw that this was true for $n=3$,
so suppose that $n\geq4$ and that $\Phi_{n-1}$ is injective. By
checking that both maps agree on the elements of $\mathscr{X}_{n-1}$, it is
easy to see that~${\Psi\circ\Phi|_{eMe}\circ\Upsilon=\Phi_{n-1}}$.
(The map $\Upsilon$ was defined at the end of Section
\ref{sect:fdsim}., and $\Psi$ in Proposition 13.) Now $\Psi$ is surjective (by
Proposition~\ref{eMegens}) and $\Phi_{n-1}$ is injective (by
assumption), so it follows that $\Phi|_{eMe}$ is injective. After
the proof of Proposition \ref{Misinverse}, we observed that
$\Phi|_{E(M)}$ and $\Phi|_{G(M)}$ are injective. By Lemma
\ref{propP}, $e$ has property (P) and it follows, by Proposition
\ref{firstprop}, that $\Phi$ is injective. \epf

We remark that the method of Propositions 5 and 13 may be used to provide a concise proof of the presentation of $\mathscr{I}_n$ originally found by Popova \cite{7}.

\end{document}